\documentclass[12pt,a4paper]{article}
\usepackage[OT1]{fontenc}
\usepackage[latin1]{inputenc}
\usepackage[english]{babel}
\usepackage{amsfonts}
\usepackage{amsthm}
\newtheorem{teo}{Theorem}[section]
\newtheorem{prop}[teo]{Proposition}
\newtheorem{lem}[teo]{Lemma}
\newtheorem{coro}[teo]{Corollary}

\theoremstyle{definition}
\newtheorem{rem}[teo]{Remark}

\def\h{{\cal H}}
\def\b{{\cal B}}

\def\be{{\cal B}_2({\cal H})}
\def\beh{{\cal B}_2({\cal H})_h}
\def\beah{{\cal B}_2({\cal H})_{ah}}
\def\gldos{Gl_2({\cal H})}
\def\u2{U_2({\cal H})}

\def\gl{Gl({\cal H})}
\def\o2{Gr_{res}^0(p)}

\begin{document}

\title{\vspace*{0cm}Hopf-Rinow Theorem in the Sato Grassmannian}
\date{}
\author{Esteban Andruchow and Gabriel Larotonda}

\maketitle

\abstract{\footnotesize{\noindent Let $U_2({\cal H})$ be the Banach-Lie group of unitary operators in the Hilbert space ${\cal H}$ which are Hilbert-Schmidt perturbations of the identity $1$. In this paper we study the geometry of the unitary orbit
$$
\{ upu^*: u\in U_2({\cal H})\},
$$
of an infinite projection $p$ in ${\cal H}$. This orbit  coincides with the connected component of $p$ in the Hilbert-Schmidt restricted Grassmannian $Gr_{res}(p)$ (also known in the literature as the Sato Grassmannian) corresponding to the polarization ${\cal H}=p({\cal H})\oplus p({\cal H})^\perp$. It is known that the components of $Gr_{res}(p)$ are differentiable manifolds. Here we give a simple proof of the fact that $Gr_{res}^0(p)$ is a smooth submanifold of the affine Hilbert space $p+{\cal B}_2({\cal H})$, where ${\cal B}_2({\cal H})$ denotes the space of Hilbert-Schmidt operators of ${\cal H}$. Also we show that $Gr_{res}^0(p)$ is a homogeneous reductive space.  We introduce a natural metric,
which consists in endowing every tangent space with the trace inner product, and consider its Levi-Civita connection. This connection has been considered before, for instance its sectional curvature has been computed. We show that the Levi-Civita connection coincides with a linear connection induced by the reductive structure, a fact which allows for the easy computation of the geodesic curves. We prove that the geodesics of the connection, which are of the
form $\gamma(t)=e^{tz}pe^{-tz}$, for $z$ a $p$-codiagonal anti-hermitic element of ${\cal B}_2({\cal H})$, have minimal length provided that $\|z\|\le \pi/2$. Note that the condition is given in terms of the usual operator norm, a fact which implies that there exist minimal geodesics of arbitrary length. Also we show that any two points $p_1,p_2\in Gr_{res}^0(p)$ are joined by a minimal geodesic. If moreover $\|p_1-p_2\|<1$, the minimal geodesic is unique. Finally, we
replace the $2$-norm by the $k$-Schatten norm ($k>2$), and prove that the geodesics are also minimal for these norms, up to a critical value of $t$, which is estimated also in terms of the usual operator norm. In the process, minimality results in the $k$-norms are also obtained for the group $U_2({\cal H})$. 
}}

\setlength{\parindent}{0cm} 

\section{Introduction}
Let $\h$ be an infinite dimensional Hilbert space and ${\cal B}(\h)$ the space
of bounded linear operators acting in $\h$. Denote by $\gl$ and $U(\h)$ the
groups of invertible and unitary operators in $\h$, and by $\be$ the space of
Hilbert-Scmidt operators, that is
$$
\be=\{ a\in {\cal B}(\h): Tr(a^*a)<\infty\},
$$
where $Tr$ is the usual trace in ${\cal B}(\h)$.
This space is a two sided ideal of ${\cal B}(\h)$, and a Hilbert space with the
inner product
$$
<a,b>=Tr(b^*a).
$$
The norm induced by this inner product is called the $2$-norm, and denoted by
$$
\|a\|_2=Tr(a^*a)^{1/2}.
$$
Throughout this paper, $\| \ \|$ denotes the usual operator norm.
Consider the following groups of operators:
$$
\gldos=\{g\in \gl: g-1\in \be\},
$$
and
$$
\u2=\{u\in U(\h): u-1\in \be\},
$$
here $1\in {\cal B}(\h)$ denotes the identity operator. These groups are
examples of what in the literature \cite{harpe} is called a {\it classical}
Banach-Lie group. They have differentiable structure when endowed with the
metric $\|g_1-g_2\|_2$ (note that $g_1-g_2\in\be$).

Fix a selfadjoint infinite projection $p\in{\cal B}(\h)$. The aim of this paper is the
geometric study of the set
$$
\o2=\{ upu^*: u\in\u2\},
$$
 the connected component of $p$ in the (Hilbert-Schmidt) restricted Grassmannian corresponding
to the polarization $\h=R(p)\oplus R(p)^\perp$ \cite{segalwilson}.
Since both $R(p)$ and
$R(p)^\perp$ are infinite dimensional, the $\u2$-orbit of $p$ lies inside $Gr_{res}^0(p)$ ,
\cite{pressleysegal}. The  fact that the group $\u2$ acts transitively on $\o2$ was
proved by Stratila and Voiculescu in \cite{stravoic} (see also  \cite{carey}).
 
The Hopf-Rinow Theorem is not valid in infinite dimensional complete manifolds: two
points in a Hilbert-Riemann complete manifold may not be joined by a minimal
geodesic \cite{grossman}, \cite{mcalpin}, or even a geodesic \cite{atkin}. The main
results in this paper establish the validity of the Hopf Rinow Theorem for $\o2$
(\ref{hr1}, \ref{hr2}). In the process we prove minimality results for $\u2$, which
are perhaps well known but for which we could find no references. We also prove
minimality results, both in $\u2$ and $\o2$, for the Finsler metrics given by the
Schatten $k$-norms ($k\ge 2$).

If one chooses unitaries $\omega_n$,
$n\in{\mathbb Z}$, in the different components of ${\cal U}_{res}(\h)$, then the
connected components of $Gr_{res}(p)$ are
$$
Gr_{res}^n(p)=\omega_n\o2 \omega_n^*=\{\omega_n upu^*\omega_n^*: u\in\u2 \},
$$
so that the results described above are valid also in the other components  of $Gr_{res}(p)$.

The restricted Grasamannian is related to several areas of mathematics and physics:
loop groops \cite{segalwilson}, \cite{pressleysegal}, integrable systems
\cite{sato}, \cite{segalwilson}, \cite{mulase}, \cite{zelikin}, group representation
theory \cite{stravoic}, \cite{carey}, \cite{segal}, string theory
\cite{alvarezgomez}, \cite{bowickrajeev}.

Unitary orbits of operators, and in particular projections, have been studied before in
(\cite{cprprojections}, \cite{pr}, \cite{belrattum}, \cite{bona}, \cite{larotonda},
\cite{andsto}). In this particular framework, restricting the action to
classical groups, certain results can be found in  \cite{carey}, \cite{belrattum} and
\cite{bona}. In the latter paper, the author considered the orbit of a finite
rank projection. 

If $p$ has infinite rank and corank, then $p$ and $\be$ are linearly
independent. We shall denote by
$$
p+\be=\{p+a: a\in \be\}.
$$
Note that every element $x$ in $p+\be$ has a unique decomposition $x=p+a$, $a\in
\be$. We shall endow this affine space with the metric induced by the $2$-norm:
if $x=p+a$ and $y=p+b$, $\|x-y\|_2:=\|a-b\|_2$. Apparently, $p+\be$ is a Hilbert
space.

The orbit $\o2$ sits inside $p+\be$:
\begin{eqnarray}
q & = & upu^*=(1+(u-1))p(1-(u-1)^*)\nonumber\\
& = &p+(u-1)p+p(u-1)^*+(u-1)p(u-1)^*\in p+\be.\nonumber
\end{eqnarray}

Therefore we shall consider $\o2$ with the topology induced by the $2$-metric of
$p+\be_2$.
Throughout this paper, $L_2$ denotes the length functional for piecewise smooth
curves, either in $\u2$ or $\o2$, measured with the $2$-norm:
$$
L_2(\alpha)=\int_{t_0}^{t_1} \|\dot{\alpha}(t)\|_2\, d t.
$$
We use the subscript $h$ (resp. $ah$) to denote the sets of hermitic (resp.
anti-hermitic.) operators, e.g. $\beah=\{x\in\be: x^*=-x\}$.

Let us describe the contents and main results of the paper.

In section 2 we prove (Theorem 2.4) that $\o2$ is a smooth submanifold of the
affine Hilbert space $p+\be$, and that the map $\pi_p:\u2\to \o2$,
$\pi_p(u)=upu^*$ is a submersion.

In section 3 we introduce a linear connection, which is the Levi-Civita
connection of the trace inner product in $\o2$. This connection was considered
in \cite{pekonen}, where its sectional curvature was computed, and shown to be
non negative. It is presented here as the connection obtained from the reductive
structure for the action of $\u2$:
$$
\beah=\{y\in\beah: py=yp\} \oplus \{z\in\beah: pzp=(1-p)z(1-p)=0\},
$$
regarded as the decomposition of the Lie algebra of $\u2$ (equal to $\beah$),
the first subspace is the Lie algebra of the isotropy group, and the second
subspace is its orthogonal complement with respect to the trace. Therefore the geodesics can be
explicitely computed.

In section 4 we prove the main  results (Theorems \ref{inicial}, \ref{hr1} and
\ref{hr2}) on minimality of
geodesics with given initial (respectively boundary) data. These
results show that any pair of points in $\o2$ can be joined by a minimal geodesic, and that there are mimimal geodesics  which have arbitrary
length.

In section 5 we consider the minimality problem, when the length of a curve is 
measured with (the Finsler metric given by) the Schatten $k$-norms, for
$2<k<\infty$. Here we obtain (Theorem 5.5) that for these metrics, the
geodesics of the connection have minimal length up to a critical value of $t$
(which depends on the usual operator norm of the initial
data). In both settings, $k=2$ and $k>2$, the minimality results are proved
first in $\u2$, and then derived for $\o2$ via a natural inmersion of
projections as unitaries (more specifically, symmetries).

\section{Differentiable structure of $\o2$}

As said above, it is known that $\o2$, being a connected component of
$Gr_{res}(p)$ (the component of virtual dimension $0$ \cite{segalwilson}), is a differentiable
manifold.
Here we show that $\o2\subset p+\be_{h}$ is a differentiable (real analytic)
submanifold.
The action of $\u2$ induces the map
$$
\pi_p:\u2\to \o2 \ , \ \ \pi_p(u)=upu^*.
$$
This map, regarded as a map on $p+\be$, is real analytic. Its differential at
the identity is the linear map
$$
\delta_p: \beah \to \beh ,\  \delta_p(x)=xp-px.
$$
Here we have identified the Banach-Lie algebra of $\u2$ with the space $\beah$ of
anti-hermitic elements in $\be$, and used the fact that $\pi_p$ takes hermitic
values, i.e. in the set $\beh$ of hermitic elements of $\be$.

\begin{lem}\label{proyecciondelta}
The map $\delta_p^2$ defines an idempotent operator acting on $\beh$. Moreover,
it is symmetric for the trace inner product in $\beh$.
\end{lem}
\begin{proof}
Straightforward computations show that if one regards $\delta_p$ as a linear map
from $\b(\h)$ to $\b(\h)$, then it verifies $\delta_p^3=\delta_p$. Therefore
$\delta_p^2$ is an idempotent, whose  range and kernel coincide with the range
and kernel of $\delta_p$. Note that
$$
\delta_p^2(x)=xp-2pxp+px.
$$
Clearly $\delta_p^2$ maps $\beh$ into $\beh$ and $\beah$ into $\beah$, so that
in particular it defines an idempotent (real) linear operator acting in $\beh$.
Finally, if $x,y\in \beh$,
\begin{eqnarray}
<\delta_p^2(x),y> &=&Tr(y(xp-2pxp+px)=Tr(pyx)-2Tr(pypx)+Tr(ypx)\nonumber\\
&=&Tr((py-2pyp+yp)x)=<x,\delta_p^2(y)>.\nonumber
\end{eqnarray}
\end{proof}
Next let us show that the map $\pi_p$ is a fibration:
\begin{prop}
The map
$$
\pi_p:\u2\to \o2\subset p+\beh
$$
has continuous local cross sections. In particular, it is a locally trivial
fibre bundle.
\end{prop}
\begin{proof}
It is well known that if $p,q$ are projections such that $\|p-q\|<1$, then
the element
$s=qp+(1-q)(1-p)$ is invertible. If  $q\in \o2$, then $s\in \gldos$. Indeed,
there exists $u\in \u2$ (i.e. a unitary such that $u_0=u-1\in\be$) such that
$q=upu^*$.
Then
\begin{eqnarray}
s & =& p+u_0p+pu_0^*p+u_0pu_0^*p +1-p +u_0(1-p)+(1-p)u_0^*(1-p)\nonumber\\
&& +u_0(1-p)u_0^*(1-p)\in1+\be.\nonumber
\end{eqnarray}
Morever, $s$ verifies $sp=qp=qs$. Let $s=w|s|$ be the polar decomposition of
$s$. Note that $sp=qs$ implies that $ps^*=s^*q$, and then $s^*s$ commutes with
$p$. Therefore
$$
wpw^*=s|s|^{-1}p|s|^{-1}s^*=s(s^*s)^{-1/2}p(s^*s)^{-1/2}s^*=sp(s^*s)^{-1}s^*=sps^{-1}=q.
$$
We claim that $w\in\u2$. Indeed, $\mathbb{C} 1+\be$ is a *-Banach algebra (it is the
unitization of $\be$) with the $2$-norm: $\|\lambda
1+a\|_2^2=|\lambda|^2+\|a\|_2^2$. Since $s$ lies in $\gldos$, in particular it
is an invertible element of this Banach algebra, and therefore by the Riesz functional 
calculus, $w=s|s|^{-1}\in \mathbb{C}+\be$, so that $w=\mu 1+w_0$ with $w_0\in\be$.
On the other hand, note that $s^*s$, it is a positive operator which lies
in the C$^*$-algebra $\mathbb{C} 1+{\cal K}(\h)$, the unitization of the ideal
${\cal K}(\h)$ of compact operators. Therefore its square root is of the form $r
1+b$, with $r\ge 0$ and $b$ compact. Then $s^*s=(r 1+b)^2=r^2 1+b'$, with $b'\in
{\cal K}(\h)$. One has that $s^*s\in\gldos$, so that it is of the form $1$ plus
a compact operator, and since $\mathbb{C} 1$ and ${\cal K}(\h)$ are linearly
independent, it follows that $r=1$.  Then $w=s|s|^{-1}$ is of the form $1$ plus
compact. By the same argument as above, this implies that $\mu=1$.

The map sending an arbitrary invertible element $g\in \gldos$ to its unitary
part $u\in\u2$ is a continuous map between these groups. In fact, it is real
analytic, by the regularity properties of the Riesz functional calculus.

Summarizing, consider the open set $\{q\in \o2: \|q-p\|_2<1\}$ in $\o2$. If
$q$ lies in this open set, then in particular $\|q-p\|\le \|q-p\|_2<1$, so
that $s$ defined above lies in $\gldos$, and therefore its unitary part
$u\in\u2$ verifies $upu^*=q$. Denote by $u=\sigma_p(q)$. Clearly $\sigma_p$ is
continuous, being the composition of continuous maps. Thus it is a continuous
local cross section for $\pi_p$ on a neighbourhood of $p$. One obtains cross
sections defined on neighbourhoods of other points of $\o2$ by translating this
one with the action of $\u2$ in a standard fashion.
\end{proof}

Let us transcribe the following result contained in the appendix of the paper
\cite{rae} by I. Raeburn, which is a consequence of the implicit function
theorem in Banach spaces.

\begin{lem}
Let $G$ be a Banach-Lie group acting smoothly on a Banach space $X$. For a fixed
$x_0\in X$, denote by $\pi_{x_0}:G\to X$ the smooth map $\pi_{x_0}(g)=g\cdot
x_0$. Suppose that
\begin{enumerate}
\item
$\pi_{x_0}$ is an open mapping, when regarded as a map from $G$ onto the orbit
$\{g\cdot x_0: g\in G\}$ of $x_0$ (with the relative topology of $X$).
\item
The differential $d(\pi_{x_0})_1:(TG)_1\to X$ splits: its kernel and range are
closed complemented subspaces.
\end{enumerate}
Then the orbit $\{g\cdot x_0: g\in G\}$ is a smooth submanifold of  $X$, and the
map
$\pi_{x_0}:G\to \{g\cdot x_0: g\in G\}$ is a smooth submersion.
\end{lem}
This lemma applies  to our situation:
\begin{teo}
The orbit $\o2$ is a real analytic submanifold of $p+\be$ and the map
$$
\pi_p: \u2\to \o2 \ , \ \ \pi_p(u)=upu^*
$$
is a real analytic submersion.
\end{teo}
\begin{proof}
In our case, $G=\u2$, $X=p+\be$, $x_0=p$. The above proposition implies that
$\pi_p$ is open. The differential $d(\pi_p)_1$ is $\delta_p$, its kernel is
complemented because it is a closed subspace of the real Hilbert space
$\beah=(T\u2)_1$. The range of $\delta_p$ equals $\delta_p^2(\beh)$, and therefore it is closed and complemented in $\be$, by Lemma \ref{proyecciondelta}. In
our context, smooth means real analytic (the group and the action are real
analytic).
\end{proof}
\section{Linear connections in $\o2$ and $\u2$}

The tangent space of $\o2 $ at $q$ is
$$
(T\o2)_q=\{xq-qx: x\in\beah\},
$$
or equivalently, the range of $\delta_q$, the differential of $\pi_p$ at $q$. As noted
above, it is a closed linear subspace of $\beh$, and the operator $\delta_q^2$
is the orthogonal projection onto $(T\o2)_q$. It is natural to consider the
Hilbert-Riemann metric in $\o2$ which consists of endowing each tangent space
with the trace inner product. Therefore the Levi-Civita connection of this
metric is given by differentiating in the ambient space $\beh$ and projecting
onto $T\o2$. That is, if $X$ is a tangent vector field along a curve $\gamma$
in  $\o2$, then
$$
\frac{D X}{d t}=\delta_\gamma^2(\dot{X}).
$$
This same connection can be obtained by other means, it is the connection
induced by the action of $\u2$ on $\o2$ and the decomposition of the Banach-Lie algebra $\beah$ of
$\u2$:
$$
\beah=\{y\in \beah: yp=py\}\oplus \{z\in \beah: pzp=(1-p)z(1-p)=0\},
$$
or, if one regards operators as $2\times 2$ matrices in terms of $p$, the
decomposition of $\beah$ in diagonal plus codiagonal matrices. This type of
decomposition, where the first subspace is the Lie algebra of the isotropy group
of the action (at $p$), and the second subspace is invariant under the inner
action of the isotropy group, is what in differential geometry is called a
reductive structure of the homogeneous space \cite{sharpe}. We do not perform
this construction here, it can be read in \cite{cprprojections}, where it is
done in a different context, but with computations that are formally identical.
This alternative description of the Levi-Civita connection of $\o2$ allows for
the easy computation of the geodesics curves of the connection. The unique
geodesic $\delta$ of $\o2$ satisfying
$$
\delta(0)=q \  \hbox{ and } \  \dot{\delta}(0)=xq-qx
$$
is given by
$$
\delta(t)=e^{tz}qe^{-tz}
$$
where $z$ is the unique codiagonal element in $\beah$ ($qzq=(1-q)z(1-q)=0$) such
that
$$
zq-qz=xq-qx.
$$
Equivalently, $z=\delta_q(xp-px)$ is the projection of $x$ in the decomposition
$$
x=y+z\in \{y\in \beah: yp=py\}\oplus \{z\in \beah: pzp=(1-p)z(1-p)=0\}.
$$

Although our main interest in this paper are projections, it will be useful to
take a brief look at the natural Riemannian geometry of the group $\u2$. Namely,
the metric given by considering real part of the trace inner product, and
therefore, the $2$-norm at each tangent space. The tangent spaces of $\u2$
identify with
$$
(T\u2)_u=u \beah= \beah u.
$$
As with $\o2$, the covariant derivative consists of differentiating in the
ambient space, and projecting onto $T\u2$. Geodesics of the Levi-Civita
connection are curves of the form
$$
\mu(t)=u e^{tx},
$$
for $u\in \u2$ and $x\in\beah$.
The exponential mapping of this connection is the  map
$$
exp: \beah\to \u2 , \ \ exp(x)=e^x.
$$

\begin{rem}\label{remarko}

\noindent
\begin{enumerate}

\item
The exponential map 
$$
exp: \beah\to \u2
$$
is surjective. This fact is certainly well known. Here is a simple proof. If
$u\in\u2$, then it has a spectral decomposition $u=p_0+\sum_{k\ge
1}(1+\alpha_k)p_k$, where $\alpha_k$ are the non zero eigenvalues of
$u-1\in\beah$. There exist $t_k\in{\mathbb R}$ with $|t_k|\le \pi$ such
that $e^{it_k}=1+\alpha_k$. One has the elementary estimate
$$
|t_k|^2(1-\frac{|t_k|^2}{12})\le |e^{it_k}-1|^2=|\alpha_k|^2,
$$
which implies that the sequence $(t_k)$ is square summable. Let $z=\sum_{k\ge 1}
it_kp_k$ (note that $p_k$ are finite rank pairwise orthogonal projections). Thus
$z\in \beah$ and clearly $e^z=u$.
\item
The exponential map is a bijection between the sets
$$
\beah\supset \{z\in\beah: \|z\|<\pi\}\to \{u\in\u2: \|1-u\|<2\}.
$$
Clearly if $z\in\beah$ with $\|z\|<\pi$, then $e^z\in\u2$ and
$\|e^z-1\|<2$.
Suppose that $u\in\u2$ with $\|u-1\|<2$. Then there exist $x\in{\cal B}(\h)$,
$x^*=-x$ and $\|x\|<\pi$, and $z\in\beah$, such that $e^x=e^z=u$. Since
$\|x\|<\pi$, $x$ equals a power series in $u=e^z$, which implies that $z$
commutes with $x$. Then $e^{z-x}=1$ and thus (note that $z-x$ is anti-hermitic)
$z-x=\sum_{k\ge 1}2k\pi i p_k$ for certain projections $p_k$. Let $z=\sum_{j\ge 1}
\lambda_j q_j$ be the spectral decomposition of $z$. Note that  $e^z=\sum_{j\ge 1}
e^{\lambda_j} q_j$, and since $x$ commutes with $e^z$, this implies that $x$ commutes
with  $q_j$, and also with $z$. Also it is clear that $q_j$ and $p_k$ also commute,
and that $q_j$ have finite ranks. Then
$$
x=\sum_{j\ge 1} \lambda_j q_j+\sum_{|k|\ge 1}2k\pi i p_k.
$$
The fact that $\|x\|<\pi$ implies that the terms $2k\pi i
p_k$ are  cancelled by some of the $\lambda_j q_j$, in order that none of the
remaining $\lambda_j$ verify $|\lambda_j|\ge \pi$. It follows that the $p_k$
have finite ranks, and that there are finitely many. Thus we can define $z'$
adding the remaining $\lambda_j q_j$. Clearly $z'$ verifies $\|z'\|<\pi$ and
$e^{z'}=e^z=u$.
\item
The  argument above in fact shows that when one considers the exponential 
$$
exp:\{z\in \beah: \|z\|\le \pi\}\to \u2,
$$
then it is surjective.

\end{enumerate}
\end{rem}

If $x\in{\cal B}(\h)$ with $\|x\|<\pi$ then it is well known that
\begin{equation}\label{arcoseno}
\|e^x-1\|=2\sin(\frac{\|x\|}{2}).
\end{equation}

There is a natural way to imbedd projections in the unitary group by means of
the map $q\mapsto \epsilon_q=2q-1$. The unitary $\epsilon_q=2q-1$ is a symmetry,
i.e. a selfadjoint unitary: $\epsilon^*=\epsilon$, $\epsilon^2=1$. See for
instance \cite{pr}, where this simple trick was used to characterize the
minimality of geodesics in the Grassmann manifold of a C$^*$-algebra. However,
if $q\in \o2$, $\epsilon_q$ does not belong to $\u2$ (recall that $q$ has
infinite rank and corank). One can slightly modify this imbedding in order that
it takes values in $\u2$. Consider:
$$
S:\o2\to \u2, \ \ S(q)=(2q-1)(2p-1)=\epsilon_q\epsilon_p.
$$
Clearly it takes unitary values. Let us show that these unitaries belong to
$\u2$. Note that if $q=upu^*$ with $u\in\u2$, then
$$
\epsilon_q=u\epsilon_pu^*=\epsilon_p+(u-1)\epsilon_p+\epsilon_p(u^*-1)+(u-1)\epsilon_p(u^*-1)\in
\epsilon_p+\be,
$$
so that
$$
\epsilon_q\epsilon_p\in (\epsilon_p+\be)\epsilon_p=1+\be.
$$
\begin{prop}\label{mapaS}
The map $S$ preserves geodesics, and its differential is $2$ times an isometry.
\end{prop}
\begin{proof}
Let $\delta$ be a geodesic in $\o2$, $\delta(t)=e^{tz}qe^{-tz}$ with $z\in\beah$
and $z$ codiagonal in terms of $q$. This latter condition is equivalent to $z$
anticommuting with $\epsilon_q$: $z\epsilon_q=-\epsilon_q z$. Which implies, as
remarked in \cite{cprprojections}, that $\epsilon_qe^{-tz}=e^{tz}\epsilon_q$,
and thus $e^{tz}\epsilon_q e^{-tz}=e^{2tz}\epsilon_q$.
Therefore
$$
S(\delta(t))=e^{2tz}\epsilon_q\epsilon_p,
$$
which is a geodesic in $\u2$.

The differential of $S$ at $q$ is given by
$$
dS_q(v)=2v\epsilon_p, \ \ v\in (T\o2)_q.
$$
Right multiplication by a fixed unitary operator is isometric in $\be$,
therefore this map is  $2$ times an isometry.
\end{proof}

\section{Minimality of geodesics}

In this section we prove that the geodesics of the linear connection have
minimal length up to a certain critical value of $t$. This could be derived from
 the general theory of Hilbert-Riemann manifolds.
 We shall prove it here, and in the process obtain a uniform lower bound for the
geodesic
radius, i.e. the radius of normal neighbourhoods.
First we need  minimality results in the group $\u2$. These results  are perhaps
well known. We include  proofs  here for we could not find 
references for them, and they are central to our argument on $\o2$.
\begin{lem}
Suppose that $x\in \beah$ has finite spectrum and $\|x\|\le\pi$, and let
$u\in\u2$. Then the (geodesic) curve $\mu(t)=ue^{tx}$, $t\in[0,1]$, has minimal
length among all piecewise smooth curves in $\u2$ joining the same endpoints.
\end{lem}
\begin{proof}
Since the action of left multiplication by $u$ is an isometric isomorphism of
$\u2$, it suffices to consider the case $u=1$. Let
$\sigma(x)=\{\lambda_0=0,\lambda_1,\dots ,\lambda_n\}$ be the spectrum of $x$.
Then $x=\sum_{i=1}^n \lambda_i p_i$ for $p_i$ finite rank projections, and
denote by $p_0$ the projection onto the kernel of $x$.  Note that
$e^{tx}=p_0+\sum_{i=1}^n e^{t\lambda_i}p_i$. Let $r_i^2=Tr(p_i)$, $i=1,\dots, n$,
and denote by ${\cal S}_i$ the sphere in $\be$ of radius $r_i$,
$$
{\cal S}_i=\{a\in \be: Tr(a^*a)=r_i^2\},
$$
with its natural Hilbert-Riemann metric induced by the (trace) inner product in
the Hilbert space $\be$.
Consider the following smooth map
$$
\Phi: \u2\to {\cal S}_1\times \dots \times {\cal S}_n , \ \ \Phi(u)=(p_1u,\dots,
p_nu).
$$
Here the product of spheres is considered with the product metric.
Apparently $\Phi$ is well defined and smooth. Note that the curve $\Phi(\mu(t))$
is a minimal geodesic of the manifold ${\cal S}_1\times \dots \times {\cal
S}_n$. Indeed,
$$
\Phi(\mu(t))=(e^{t\lambda_1}p_1,\dots, e^{t\lambda_n}p_2),
$$
where each coordinate $e^{t\lambda_i}p_i$ is a geodesic of the corresponding
sphere ${\cal S}_i$, with length equal to $|\lambda_i|r_i\le \|x\|r_i\le\pi
r_i$, and therefore it is minimal. Then $\Phi(\mu(t))$ is minimal, being the
cartesian product of $n$ minimal geodesics in the factors.
Next, note that the length of $\Phi(\mu)$ equals the length of $\mu$:
$$
L_2(\Phi(\mu))=\int_0^1 \|(\lambda_1 e^{t\lambda_1}p_1,\dots,\lambda_n
e^{t\lambda_n}p_n)\| d t=\{\sum_{i=1}^n
|\lambda_i|^2r_i^2\}^{1/2}=\|x\|_2=L_2(\mu).
$$
If $\nu(t)$, $t\in[0,1]$ is any other smooth curve in $\u2$, we claim that
$L_2(\Phi(\nu))\le L_2(\nu)$. Clearly this would prove the lemma. Indeed, since
$\Phi(\mu)$ is minimal in ${\cal S}_1\times \dots \times {\cal S}_n$, one has
$L_2(\Phi(\mu))\le L_2(\Phi(\nu))$, and therefore
$$
L_2(\nu)\ge L_2(\Phi(\nu))\ge L_2(\Phi(\mu))=L_2(\mu).
$$
Note that
$$
L_2(\Phi(\nu))=\int_0^1\{\sum_{i=1}^n\|\dot{\nu}p_i\|_2^2\}^{1/2} dt .
$$
Since $\sum_{i=1}^np_i=1-p_0$ and $\dot{\nu}^*(1-p_0)\dot{\nu}\le
\dot{\nu}^*\dot{\nu}$, one has that
$$
\sum_{i=1}^n\|\dot{\nu}p_i\|_2^2=\sum_{i=1}^nTr(\dot{\nu}^*p_i\dot{\nu})=
Tr(\dot{\nu}^*(1-p_0)\dot{\nu})\le Tr(\dot{\nu}^*\dot{\nu})=\|\dot{\nu}\|_2^2.
$$
Therefore
$$
L_2(\Phi(\nu))\le\int_0^1 \|\dot{\nu}\|_2 d t =L_2(\nu).
$$
\end{proof}

\begin{teo}\label{minimalidadunitariainicial}
Let $u\in\u2$ and $x\in\beah$ with $\|x\|\le \pi$. Then the curve
$\mu(t)=ue^{tx}$, $t\in[0,1]$ is shorter than any other pieceise smooth curve in
$\u2$ joining the same endpoints. Moreover, if $\|x\|<\pi$, then $\mu$ is unique with
this property.
\end{teo}
\begin{proof}
Again, by the same argument as in the previous lemma, we may suppose $u=1$.
Assume that $\mu$ is not minimal. Let
$\gamma(t)$, $t\in[0,1]$ be a piecewise smooth curve in $\u2$ with
$L_2(\gamma)+\delta=\|x\|_2=L_2(\mu)$, for some $\delta>0$. Let $z\in \beah$
be a finite rank operator close enough to $x$ in the $2$-norm in order that
$$
y=log(e^{-x}e^z) \hbox{ verifies } \|y\|_2<\delta/4,
$$
$$
|\|x\|_2-\|z\|_2|<\delta/4,
$$
and
$$
\|z\|<\pi.
$$
Let $\rho(t)=e^xe^{ty}$, and consider $\gamma\#\rho$ the curve $\gamma$ followed
by $\rho$, which joins $1$ to $e^{z}$.
Then
$$
L_2(\gamma\#\rho)=L_2(\gamma)+L_2(\rho)=L_2(\gamma)+\|y\|_2<L_2(\gamma)+\delta/4=\|x\|_2-3\delta/4<\|z\|_2-\delta/2,
$$
which contradicts the minimality of the curve $e^{tz}$ proved in the previous
lemma, because $\|z\|\le\pi$.

Suppose now that $\|x\|<\pi$.
By the general theory of Hilbert-Riemann manifolds \cite{lang}, any minimal
curve starting at $u$ is a geodesic of the linear connection, i.e. a curve of
the form $ue^{tw}$. If it joins the same endpoints as $\mu$, then it must be
$e^w=e^{x}$. Since $\|x\|<\pi$, $x$ is a power series in terms of
$e^x$, and therefore $w$ commutes with $x$. Then $e^{w-x}=1$. Suppose that $w\ne
x$, then
$$
w-x=\sum_{k=1}^m 2k\pi i p_i,
$$
for certain pairwise orthogonal (non nil) projectors $p_i\in\be$.
Then,
$$
\|w-x\|_2^2=\sum_{k=1}^m 4\pi^2 Tr(p_i)^2\ge 4\pi^2.
$$
Since $\|x\|<\pi$, this inequality clearly impies that $\|w\|\ge \pi$,
therefore leading to a contradiction.
\end{proof}

\begin{rem}
The proof in  Theorem \ref{minimalidadunitariainicial}  shows  that  if
$x\in \beah$, the curve $e^{tx}$ remains minimal as long as $t\|x\|\le \pi$. One has
coincidence $\|x\|=\|x\|_2$ only for rank one operators. In general, the number
$C_x=\|x\|_2/\|x\|$ can be arbitrarily large. Therefore, for a specific
$x\in\beah$, in terms of the $2$-norm, $e^{tx}$ will remain minimal as long as
$$
t\|x\|_2\le C_x \pi.
$$
\end{rem}

\begin{coro}
There are in $\u2$ minimal geodesics of arbitrary length. Thus the Riemannian
diameter of $\u2$ is infinite.
\end{coro}

\begin{teo}\label{minimalidadunitariaborde}
Let $u_0, u_1 \in\u2$. Then there exists a minimal geodesic curve joining them.
If $\|u_0-u_1\|<2$, then this geodesic is unique.
\end{teo}
\begin{proof}
Again, using the isometric property of the left action of
$\u2$ on itself, we may suppose $u_0=1$.
The first assertion follows from the surjectivity of the exponential map
$exp:\{x\in\beah: \|x\|\le \pi\}\to\o2$ in Remark \ref{remarko}, and Theorem
\ref{minimalidadunitariainicial}. The uniqueness assertion also follows from Remark
\ref{remarko}. 
\end{proof}
Denote by $d_2$ the geodesic distance, i.e. the metric induced by the $2$-norm on
the tangent spaces, both in $\u2$ and $\o2$.

\begin{prop}
If $u,v\in\u2$ then
$$
\sqrt{ 1-\frac{\pi^2}{12} } \; d_2(u,v) \le \|u-v\|_2\le d_2(u,v).
$$
In particular the metric space $(\u2, d_2)$ is complete.
\end{prop}
\begin{proof}
Since left multiplication by $v^*$ is an isometry for both metrics, we may assume
that $v=1$. As in Remark \ref{remarko}, we may assume that $u=p_0 +\sum_{k\ge
1}e^{it_k}p_k$, with $p_i$ mutually orthogonal projections and $\mid t_k\mid \le
\pi$. Then
$$
\|u-1\|_2^2=\|\sum_{k\ge 1} (e^{it_k}-1)p_k\|_2^2=\sum_{k\ge 1}\mid
e^{it_k}-1\mid^2 r_k^2=\sum_{k\ge 1} 2(1-\cos(t_k)) r_k^2,
$$
where $r_k=Tr(p_k)$.
Now 
$$\mid t \mid ^2 \ge 2(1-\cos(t))\ge \mid t\mid^2 (1-\frac{\mid t\mid^2}{12} )\ge 
\mid t\mid^2 (1-\frac{ \pi^2}{12} )
$$
for any $t\in [-\pi,\pi]$. Let $z=\sum_{k\ge 1}i t_k p_k$; clearly $z\in \beah$ by
the inequality above and $e^z=u$. If $\gamma(t)=e^{tz}$, then $\gamma$ is a minimal
geodesic in $\u2$ joining $1$ to $u$ because $\|z\|\le \pi$. Then
$d_2(u,1)=L_2(\gamma)=\|z\|_2$, and from the two inequalities above we obtain
$\sqrt{ 1-\frac{\pi^2}{12} } \; \|z\|_2 \le \|u-v\|_2\le \|z\|_2$, which proves the
assertion of the proposition.

Therefore, $\u2$ is complete with the geodesic distance, because $(\u2,\| \ \|_2)$
is complete. This fact is certainly well known. We include a short proof. Suppose
that $u_n$ is a Cauchy sequence in $\u2$ for the $2$-norm. Since the $2$-norm bounds
the operator norm, it follows that there exists a unitary operator $u$ such that
$\|u_n-u\|\to 0$. On the other hand $u_n-1$ is a Cauchy sequence in $\be$, and
therefore it converges to some operator in $\be$. Thus $u\in\u2$.
\end{proof}

\begin{rem}
If $x,y$ are anti-hermitic operators with $\|x\|,\|y\|<\pi$, then $e^x=e^y$ implies
$x=y$. If $\|x\|=\|y\|=\pi$, from Theorem \ref{minimalidadunitariainicial}, it
follows that $e^x=e^y$ implies $\|x\|_2=\|y\|_2$, because the curves $e^{tx},
e^{ty}$ are both minimal geodesics joining the same endpoints, hence they have the
same length.
\end{rem}

Now our main results on minimal geodesics of $\o2$ follow:

\begin{teo}\label{inicial}
Let $z\in\beah$ which is codiagonal with respect to $q\in \o2$, and such that
$\|z\|\le\pi/2$. Then the geodesic $\alpha(t)=e^{tz}qe^{-tz}$, $t\in[0,1]$ has
minimal length among all piecewise smooth curves in $\o2$ joining the same
endpoints. Moreover, if $\|z\|<\pi/2$, then $\alpha$ is unique having this property.
\end{teo}
\begin{proof}
Let $\beta$ be any other piecewise smooth curve in $\o2$ having the same
endpoints as $\alpha$. Consider $S(\alpha)$ and $S(\beta)$ in $\u2$. Note that
$S(\alpha)(t)=e^{2tz}\epsilon_q\epsilon_p$, with  $2\|z\|\le \pi$. Therefore
$$
L_2(\alpha)=\frac12 L_2(S(\alpha))\le \frac12 (S(\beta))=L_2(\beta).
$$
The uniqueness part is an easy consequence.
\end{proof}
\begin{rem}
Again, as remarked after \ref{minimalidadunitariainicial}, for specific $z$ (of
rank greater than one), the geodesic $e^{tz}qe^{tz}$ will remain minimal  as
long as
$$
t\|z\|_2\le \frac{\pi}{2} C_z,
$$
where again $C_z=\frac{\|z\|_2}{\|z\|}$ can be arbitrarily large.
\end{rem}
Analogously as for $\u2$, one has
\begin{coro}
There are in $\o2$ minimal geodesics of arbitrary length, thus $\o2$ has infinite
Riemannian diameter.
\end{coro}

\begin{teo}\label{hr1}
Let $q_0,q_1\in \o2$ such that $\|q_0-q_1\|<1$.
Then there exists a unique geodesic joining them, which has minimal length.
\end{teo}
\begin{proof}
The action of $\u2$ on $\o2$ is isometric, therefore we may suppose without loss
of generality that $q_0=p$.
Then \cite{pr}, \cite{cprprojections} there exists $z\in{\cal B}(\h)$, $z^*=-z$,
$\|z\|<\pi/2$, $z$ $p$-codiagonal such that $e^zpe^{-z}=q_1$. Therefore
$$
\epsilon_{q_1}=e^z\epsilon_pe^{-z}=e^{2z}\epsilon_p,
$$
and thus $z=\frac12 log(\epsilon_{q_1}\epsilon_p)$, where $log$ is well defined
because
$\|1-\epsilon_{q_1}\epsilon_p\|=\|\epsilon_{q_1}-\epsilon_p\|=2\|p-q_1\|<2$.
On the other hand $\epsilon_{q_1}\epsilon_p\in\u2$, therefore by Remark
\ref{remarko}, $z\in\beah$. Moreover,
$$
2\|z\|\le \pi,
$$
and therefore the curve $\mu(t)=e^{2tz}\epsilon_p$ is a minimal geodesic in
$\u2$. Again, as in the previous theorem, this implies that the geodesic curve
$\delta(t)=e^{tz}pe^{-tz}$, which joins $p$ and $q_1$, is minimal.
\end{proof}
Next let us consider the case when $\|q_1-q_2\|=1$. The problem of existence of minimal curves in
this case, in the context of abstract $C^*$-algebras, and measuring with the operator
norm, has been studied by Brown in \cite{brown}.

\begin{teo}\label{hr2}
Let $q_0,q_1\in\o2$ with $\|q_0-q_1\|=1$. Then there exists a minimal geodesic
joining them.
\end{teo}
\begin{proof}
Again, without loss of generality, we may suppose $q_0=p$.
Consider the following subspaces:
$$
H_{00}=\ker p\cap \ker q_1 \ , \ \ H_{01}=\ker p\cap R(q_1) \ , \ \ H_{10}=R(p)\cap
\ker q_1\ , \ \ H_{11}=R(p)\cap R(q_1) \ ,
$$
and
$$
H_0=(H_{00}\oplus H_{01} \oplus H_{10} \oplus H_{11})^\perp.
$$
These are the usual subspaces to regard when considering the unitary equivalence of
two projections \cite{dixmier}. The space $H_0$ is usually called the generic part
of $H$. It is invariant both for $p$ and $q_1$. Also it is
clear that $H_{00}$ and $H_{11}$ are invariant for $p$ and $q_1$, and that $p$ and
$q_1$ coincide here. Thus in order to find a unitary operator $e^z$ conjugating $p$
and $q_1$, with $z\in \beah$,  which is codiagonal with respect to $p$, and such that
$\|z\|\le \pi/2$, one needs to focus on the subspaces $H_0$ and $H_{01}\oplus H_{10}$.

Let us treat first $H_0$, denote by $p'$ and $q'_1$ the projections $p$ and $q_1$ reduced to $H_0$. These projections are in what in the literature is called generic position. In \cite{halmos} Halmos showed that two projections in generic position are unitarily equivalent, more specifically, he showed that there exists a unitary operator $w:H_0\to K\times K$ such that 
$$
wp'w*=p''=\left( \begin{array}{cc} 1 & 0 \\ 0 & 0 \end{array} \right) \ \ \hbox{ and } \ \ 
wq'_1w*=q''_1=\left( \begin{array}{cc} c^2 & cs \\ cs & s^2 \end{array} \right),
$$
where $c,s$ are positive commuting contractions acting in $K$ and satisfying $c^2+s^2=1$. We claim that there exists an anti-hermitic operator $y$ acting on $K\times K$, which is a co-diagonal matrix, and  such that $e^yp''e^{-y}=q''_1$. In that case, the element $z_0=w^*yw$ is an anti-hermitic operator in $H_0$, which verifies $e^{z_0}p'e^{-z_0}=q'_1$, and is co-diagonal with respect to $p'$. Moreover. we claim that $y$ is a Hilbert-Schmidt operator in $K\times K$ with $\|y\|\le \pi/2$, so that $z_0$ is also a Hilbert-Schmidt operator in $H_0$ with $\|z_0\|\le\pi/2$. Let us prove these claims. By a functional calculus argument, there exists a positive element $x$ in the C$^*$ algebra generated by $c$, with $\|x\|\le \pi/2$,  such that $c=cos(x)$ and $s=sin(x)$. Since $q''_1$ lies in the Hilbert-Schmidt Grassmannian of $p''$, in particular one has that $q''_1|_{R(p'')}$ is a Hilbert-Schmidt operator. That is, the operator $cos(x)\ sin(x)+sin(x)^2$ is Hilbert-Schmidt in $K$. By a strightforward functional calculus argument, it follows that $x$ is a Hilbert-Schmidt operator. Consider the operator
$$
y=\left( \begin{array}{cc} 0 & -x \\ x & 0 \end{array} \right)
$$
Clearly $y^*=-y$, $\|y\|\le \pi/2$. A straightforward computation shows that 
$$
e^yp''e^{-y}=p''_1,
$$
and our claims follow.

Let us consider now the space  $H_{01}\oplus H_{10}$. Recall
\cite{segalwilson} that an alternative definition of $\o2$ states that if  $q_1\in \o2$ then 
$$
pq_1|_{R(q_1)}:R(q_1)\to R(p)
$$
is a Fredholm operator of index $0$. Note that $H_{01}=\ker(pq_1|_{R(q_1)})$. Thus
in particular $\dim H_{01}<\infty$.
On the other hand, it is also apparent that $H_{10}\subset R(pq_1)^\perp \cap R(p)$,
and therefore also $\dim H_{10}<\infty$. Therefore, the fact that $pq_1|_{R(q_1)}$
has zero index implies that 
$$
\dim H_{01}\le \dim H_{10}.
$$
The fact that $q_1$   lies in the connected component of $p$ in the Sato Grassmannian
corresponding to the polarization given by $p$, implies that, reciprocally, $p$ lies
in the component of $q_1$, in the Grassmannian corresponding to the polarization
given by $q_1$. Thus, by symmetry, 
$$
\dim H_{01}= \dim H_{10}.
$$
Let $v:H_{10}\to H_{01}$ be a surjective isometry, and consider 
$$
w: H_{01}\oplus H_{10}\to H_{01}\oplus H_{10} \ , w(\xi'+\xi'')=v^*\xi'+v\xi''.
$$
In matrix form (in terms of the decomposition $H_{01}\oplus H_{10}$), 
$$
w=
\left(
\begin{array}{ll}
0 & v \\ v^* & 0
\end{array}
\right).
$$
Apparently, $w p|_{H_{01}\oplus H_{10}} w^*=q_1|_{H_{01}\oplus H_{10}}$. Let 
$$
z_{2}= \pi/2
\left(
\begin{array}{ll}
0 & v \\ -v^* & 0
\end{array}
\right).
$$
Note that $z_2$ is an anti-hermitic operator in $H_{01}\oplus H_{10}$, with norm
equal to $\pi/2$.
 A straightforward matrix computation shows that $e^{z_2}=w$. Consider now
 $$
 z=z_0+z_1+z_2,
 $$
 where $z_1=0$ in  $H_{00}\oplus H_{11}$, and $z_0$ is the anti-hermitic operator in
the generic part $H_0$ of $H$ found above. Then it is clear that $z$ is
anti-hermitic, Hilbert-Schmidt ($\dim(H_{01}\oplus H_{10})<\infty$), $p$-codiagonal,
$\|z\|=\pi/2$, and  $e^zpe^{-z}=q_1$.
\end{proof}

Also completeness of the geodesic metric follows:
\begin{coro}
The metric space $(\o2, d_2)$ is complete.
\end{coro}
\begin{proof}
Let $q_n$ be a Cauchy sequence in $\o2$. Since the map $S:\o2\to \u2$ of Proposition
\ref{mapaS} is $2$ times an isometry, then $S(q_n)$ is a Cauchy sequence in $\u2$,
and therefore converges to an element $u$ of $\u2$ in the metric $d_2$. Moreover,
$S(\o2)$ is closed in $\u2$ and then exists $q\in \o2$ such that $S(q)=u$. Clearly
$d_2(q_n,q)\to 0$.
\end{proof}

\section{$k$-norms}
In this section we study the minimality problem of geodesics in $\o2$ measured
in
the $k$-norms, for $k\in{\mathbb R}$, $k>2$. To do this, as with the case
$k=2$, we study first short curves in $\u2$ with these norms. Minimality of
geodesics in $\o2$ will follow with
arguments
similar as in the previous section. We shall endow now the tangent spaces of
$\u2$ and $\o2$ with the Schatten $k$-norm:
$$
\|x\|_k=Tr(|x|^k)^{1/k}=Tr((x^*x)^{k/2})^{1/k}.
$$
Note that since the tangent spaces live inside $\be$, and $k>2$,  the
$k$-norm of $x$ is finite. We shall denote by $L_k$ the functional which measures
the length of a curve (either in $\u2$ or $\o2$) in the $k$-norm:
$$
L_k(\alpha)=\int_{t_o}^{t_1} \|\dot{\alpha}(t)\|_k d t .
$$
We are now then in the realm of (infinite dimensional) Finsler geometry.

To prove our results, two inequalities proved by Hansen and Pedersen in
\cite{pedersenhansen} will play a fundamental role. Let us transcribe these
inequalities, called Jensen's  inequalities.

\noindent

\begin{enumerate}
\item
The first is the version  for C$^*$-algebras  (\cite{pedersenhansen}, Th. 2.7):
if  $f(t)$ is a convex continuous
real  function, defined on an interval $I$ and and $A$ is a C$^*$-algebra with
finite unital trace $tr$, then the inequality
\begin{equation}
tr\bigl( f(\sum_{i=1}^n b_i^*a_ib_i)\bigr)\le tr\bigl( \sum_{i=1}^n
b_i^*f(a_i)b_i\bigr)
\end{equation}
is valid for every $n$-tuple $(a_1,\dots,a_n)$ of selfadjoint elements in $A$
with spectra contained in $I$ and every $n$-tuple $(b_1,\dots,b_n)$ in $A$ with
$\sum_{i=1}^n b_i^*b_i=1$.
We shall use it in a simpler form:
if $a$ is a selfadjoint  element in a $C^*$-algebra with  trace $tr$, then
\begin{equation}\label{jensen}
tr(f(a))\le f(tr(a))
\end{equation}
for every convex continuous real function defined in the
spectrum of $a$.
\item
The second inequality is valid for finite matrices (\cite{pedersenhansen}, Th.
2.4): let $f$ be a convex continuous function defined on $I$ and let $m$ and $n$
be natural numbers, then
\begin{equation}\label{jensenmatrices}
Tr(f(\sum_{i=1}^n a_i^*x_ia_i))\le Tr(\sum_{i=1}^n a_i^*f(x_i) a_i)
\end{equation}
for every $n$-tuple $(x_1,\dots,x_n)$ of selfadjoint $m\times m$ matrices with
spectra contained in $I$ and every $n$-tuple $(a_1,\dots,a_n)$ of $m\times m$
matrices with $\sum_{i=1}^n a_i^*a_i=1$.
We shall need a simpler version, namely if $r\in{\mathbb R}$,   $r\ge 1$,
 then
 \begin{equation}\label{jensenmatricessimplificada}
Tr(a^r)=Tr( (\sum_{j=0}^n p_ja^rp_j)\ge Tr( (\sum_{j=0}^n
p_jap_j)^r)=Tr(\sum_{j=0}^n (p_jap_j)^r),
\end{equation}
for $p_0,p_1,\dots,p_n$ projections with $\sum_{j=0}^n p_j=1$ and $p_1,\dots,
p_n$ of finite rank, and $a$ a positive trace class operator. A simple
aproximation argument shows that one can obtain (\ref{jensenmatricessimplificada})
from (\ref{jensenmatrices}). Indeed, let $\{\xi_1^j,\dots,\xi_{k_j}^j\}$ be an
orthonormal basis for the range of $p_j$, $j=1,\dots, n$ and
$\{\varphi_i,\varphi_2,\dots \}$ be an orthonormal basis for the range of 
$p_0$. For any integer $N\ge 1$, let $e_N$ denote the orthogonal projection onto
the subspace generated by $\{\xi_i^j, j=1,\dots, n, \ i=1,\dots
k_j\}\cup\{\varphi_k, k=1,\dots, N\}$. Clearly $e_N$ is a finite rank projection
such that $p_j\le e_N$ , for $j=1,\dots,n$ and such that $e_Np_0e_N=p_{0,N}$ is
also a projection. Let $a_N=e_Nae_N$. Then the following facts are apparent:
\begin{enumerate}
\item
$p_ja_Np_j=p_jap_j$ for $j=1,\dots, n$.
\item
$a_N\to a$ and $p_{0,N}a_Np_{0,N}\to p_0ap_0$ in  $\| \ \|_1$, and therefore
$p_ja_N^rp_j\to p_ja^rp_j$ for $j=0,1,\dots,n$ and $r\ge 1$ in $\| \ \|_1$.
\end{enumerate}
It follows that one can reduce to prove (\ref{jensenmatricessimplificada}) for the
operator $a_N$ and the projections $p_{0,N},p_1,\dots, p_n$, all of  which are operators
in the range of $e_N$, which is finite dimensional.

\end{enumerate}
Let us first state  the following  lemma which is a simple consequence of
(\ref{jensen}).
\begin{lem}\label{lemaresucitado}
Let $a\in \b(\h)$ be a positive  operator  and $p$ a finite rank projection.
Then, if $r\in{\mathbb R}$, $r\ge 1$
$$
Tr(pap)^r\le Tr(p)^{r-1}Tr((pap)^r).
$$
\end{lem}
\begin{proof}
If $p=0$ the result is trivial. Suppose $Tr(p)\ne 0$.
Consider the finite C$^*$-algebra $p\b(\h) p$, with unit $p$ and normalized
finite trace
$tr(pxp)=\frac{Tr(pxp)}{Tr(p)}$.
Then by Jensen's trace inequality for the map $f(t)=t^r$,
$$
\frac{Tr(pap)^r}{Tr(p)^r}\le \frac{Tr((pap)^r)}{Tr(p)},
$$
which is the desired inequality.
\end{proof}
Denote by ${\cal S}_R^k$ the unit sphere of $\b_k(\h)$:
$$
{\cal S}_R^k=\{x\in\b_k(\h): \|x\|_k=R\}.
$$
If $\mu(t)$ is a curve of unitaries in $\u2$, and $p$ is finite rank projection
with $Tr(p)=R^k$, then $\mu(t)p$ is a curve in ${\cal S}_R^k$:  $\|\mu
p\|_k=Tr((p\mu^*\mu p)^{k/2})^{1/k}=R$.

\begin{lem}
Let $p$ be a finite rank projection with $Tr(p)=R^k$ and $\mu(t)$ be a smooth
curve  in $\u2$, such that $\mu(0)p=p$ and
$\mu(1)p=e^{\alpha}p$ with $-\pi\le\alpha \le \pi$. Then the curve  $\mu p$
of
${\cal S}_R^k$, measured with the $k$-norm, is longer than the curve
$\epsilon(t)=e^{it\alpha}p$.
\end{lem}
\begin{proof}
The length of $\mu p$ is (in the $k$-norm) measured by
$$
\int_0^1 \|\dot{\mu}(t) p\|_k dt=\int_0^1
Tr((p\dot{\mu}(t)^*\dot{\mu}(t)p)^{k/2})^{1/k} dt.
$$
by the inequality in the above lemma,
$$
L_k(\mu p)\ge Tr(p)^{\frac{1-k/2}{k}} \int_0^1
Tr(p\dot{\mu}(t)^*\dot{\mu}(t)p)^{1/2}dt.
$$
This last integral measures the length of the curve $\mu p$ in the $2$-sphere
${\cal S}_{R^{k/2}}^2$ of radius $R^{k/2}$ in the Hilbert space $\be$. The curves
$\epsilon(t)=e^{t\alpha}p$ are minimizing geodesics of these
spheres,
provided that $|\alpha| R^{k/2}\le\pi R^{k/2}$, which holds because
$|\alpha|\le\pi$. It follows that
$$
\int_0^1  Tr(p\dot{\mu}(t)^*\dot{\mu}(t)p)^{1/2}\ge
L_2(\epsilon)=|\alpha|Tr(p)^{1/2}.
$$
Then
$$
L_k(\mu p)\ge |\alpha|Tr(p)^{\frac{1-k/2}{k}}
Tr(p)^{1/2}=|\alpha|R=L_k(\epsilon).
$$
\end{proof}

\begin{lem}
Let $x\in\beah$ with finite spectrum, $x=\sum_{i=1}^n\alpha_i p_i$ with 
$\sum_{i=1}^n p_i=1-p_0$ ($p_0$ the kernel projection of $x$)  and
$-\pi \le\alpha_i\le\pi$ (i.e. $\|x\|\le\pi$). Then the curve
 $\delta(t)=e^{itx}$, $t\in[0,1]$ is the shortest curve in $\u2$ joining its
endpoints, when measured with the $k$-norm.
\end{lem}
\begin{proof}
Let $Tr(p_i)=R_i^k$, $i=1,\dots,n$. Note that the kernel projection $p_0$ has
infinite rank.
The length $L_k(\mu)$ of $\mu$ is measured by
$\int_0^1 \|\dot{\mu}(t)\|_k dt$.
Then, by inequality (\ref{jensenmatricessimplificada}), with
$a=\dot{\mu}(t)^*\dot{\mu}(t)\ge 0$ in $\b_1(\h)$, one has
\begin{equation}\label{desigualdad}
\|\dot{\mu}(t)\|_k\ge \{\sum_{j=0}^n
Tr\bigl((p_i\dot{\mu}(t)^*\dot{\mu}(t)p_i)^{k/2}\bigr)\}^{1/k}=
\{\sum_{i=0}^n\|\dot{\mu}(t)p_i\|_k^k\}^{1/k}.
\end{equation}

On the other hand, note that
$$
\|\dot{\delta}\|_k=\{\sum_{i=1}^n
|\alpha_i|^k R_i^k)\}^{1/k}.
$$
Trivially, $\{\sum_{j=0}^n
Tr\bigl((p_i\dot{\mu}(t)^*\dot{\mu}(t)p_i)^{k/2}\bigr)\}^{1/k}\ge \{\sum_{j=1}^n
Tr\bigl((p_i\dot{\mu}(t)^*\dot{\mu}(t)p_i)^{k/2}\bigr)\}^{1/k}$, (i.e. we omit
the term corresponding to the projection $p_0$, which has infinite trace).
We finish the proof by establishing that
$$
\{\sum_{j=1}^n
Tr\bigl((p_i\dot{\mu}(t)^*\dot{\mu}(t)p_i)^{k/2}\bigr)\}^{1/k}
\ge \{\sum_{i=1}^n
|\alpha_i|^k R_i^k)\}^{1/k}=L_k(\delta).
$$
There is a classic Minkowski type inequality (see
inequality {\bf 201} of \cite{hardylittlewoodpolya}) which
states that if $f_1,\dots,f_n$ are non negative functions, then
$$
\int_0^1 \{\sum_{i=1}^n f_i^k(t)\}^{1/k} dt \ge \bigl(\sum_{i=1}^n\{\int_0^1
f_i(t)\}^k\bigr)^{1/k}.
$$
In our case $f_i(t)=\|\dot{\mu}(t)p_i\|_k$:
$$
\int_0^1 \{\sum_{i=1}^n \|\dot{\mu}(t)p_i\|_k^k)\}^{1/k} dt \ge
\bigl(\sum_{i=1}^n\{\int_0^1 \|\dot{\mu}(t)p_i\|_k dt\}^k\bigr)^{1/k}\ge
\{\sum_{i=1}^n |\alpha_i|^k R_i^k\}^{1/k},
$$
where in the last inequality we use the previous lemma: $\int_0^1
\|\dot{\mu}(t)\|_k dt\ge |\alpha_i| R_i$ for
$i=1,\dots ,n$.
\end{proof}

\begin{teo}\label{minimalidadk}
Let $x\in \beah$  with $\|x\|\le\pi$, and $v\in \u2$. Then
the curve
$\delta(t)=ve^{tx}$  has minimal length among piecewise smooth curves in $\u2$
 joining the same endpoints, measured with the $k$-norm.
\end{teo}
\begin{proof}
There is no loss of generality if we suppose $v=1$. Indeed, for any curve $\mu$
of unitaries, $L_k(\mu)=L_k(v^*\mu)$.
Suppose that there exists a piecewise $C^1$ curve of unitaries $\mu$ which is
strictly shorter than $\delta$,
$L_k(\mu)<L_k(\delta)-\epsilon=\|x\|_k-\epsilon$. The element $x$ can be
approximated in the $k$-norm topology of $\b_k(\h)$ by anti-hermitic elements 
$z\in\b_k(\h)$, with finite spectrum and the following conditions:
\begin{enumerate}
\item
$\|z\|\le \|x\|\le\pi$.
\item
$\|x\|_k- \epsilon/2<\|z\|_k\le \|x\|_k$.
\item
There exists a $C^\infty$  curve of unitaries joining $e^{x}$ and $e^{z}$ of
$k$-length $L_k$ less than $\epsilon/2$.
\end{enumerate}
The first two are clear. The third condition can be obtained as follows. By
the third condition $e^{-x}e^{z}=e^{y}$, with $y\in \beah$. Moreover $z$ can
be adjusted so as to obtain $y$ of arbitrarily small $k$-norm. Then the curve of
unitaries $\gamma(t)=e^{x}e^{ty}$ is $C^\infty$,  joins $e^{x}$ and $e^{z}$,
with  $k$-length $\|y\|_k<\epsilon/2$.

Consider now the curve $\mu'$, which is the curve $\mu$ followed by the curve
$e^{x}e^{ty}$  above. Then clearly
$$
L_k(\mu')\le L_k(\mu)+\|y\|_k<L_k(\mu)+\epsilon/2 .
$$
Therefore $L_k(\mu')< \|x\|_k-\epsilon/2$. On the other hand, since $\mu'$
joins $1$ and $e^{z}$, by the lemma above, it must have length greater than or
equal to $\|z\|_k$. It follows that
$$
\|z\|_k\le \|x\|_k-\epsilon/2 ,
$$
a contradiction.

\end{proof}

One obtains minimality of geodesics in $\o2$ for the $k$-norm analogously as in
the previous section:
\begin{teo}
Let $z\in\beah$, codiagonal with respect to $q\in\o2$, with $\|z\|\le \pi/2$. Then
the geodesic
$\alpha(t)=e^{tz}qe^{-tz}$, $t\in[0,1]$, has minimal length for the $k$-norm 
among all piecewise smooth curves in $\o2$ having the same endpoints. If
$\|z\|<\pi/2$, this curve $\alpha$ is unique with this property.
\end{teo}
\begin{proof}
The proof follows as in the analogous result for the $2$ norm in the previous
section, noting that the map $S$ is also isometric for the $k$-norms.
\end{proof}
\begin{teo}
Let $q_1,q_2\in\o2$, then there exists a  geodesic joining them, which has minimal
length fot the $k$-norm.
\end{teo}
\begin{proof} The proof follows as in the above result, the geodesic
$\alpha(t)=e^{tz}q_1e^{-tz}$ with $\|z\|\le 
\pi/2$ exists by virtue of (\ref{hr1}) and (\ref{hr2}).
\end{proof}

\bigskip

\noindent
Esteban Andruchow and Gabriel Larotonda\\
Instituto de Ciencias \\
Universidad Nacional de Gral. Sarmiento \\
J. M. Gutierrez 1150 \\
(1613) Los Polvorines \\
Argentina  \\
e-mails: eandruch@ungs.edu.ar, glaroton@ungs.edu.ar

\end{document}